\documentclass[12pt]{article}

\usepackage[T1]{fontenc}
\usepackage[active]{srcltx}
\usepackage{lmodern}
\usepackage{amsmath,amstext,amssymb,graphicx,graphics,color}
\usepackage{mathtools}
\usepackage{enumitem}
\usepackage{tikz}
\usepackage{pgfplots}
\usepackage{theorem}
\usepackage{cite}               
\usepackage{hyperref}           
\usepackage{bm}
\usepackage{bbm}                
\usepackage{fancybox}           
\usepackage[section]{placeins}  
\usepackage{subcaption}
\theorembodyfont{\normalfont\slshape} 
\newtheorem{definition}{Definition}
\newtheorem{theorem}{Theorem}
\newtheorem{proposition}{Proposition}[section]

\newtheorem{lemma}[proposition]{Lemma}
\theoremstyle{break} 

{{\par\noindent \bf Proof. \nobreak}}%
{\nobreak \removelastskip \nobreak \hfill $\Box$ \medbreak}
{{\par\noindent \bf Proof \nobreak}}%
{\nobreak \removelastskip \nobreak \hfill $\Box$ \medbreak}
{{\par\noindent \bf Proof lemma. \nobreak}}%
{\nobreak \removelastskip \nobreak \bf End proof lemma. \medbreak}
\newenvironment{remark}{\par \medskip \noindent {\bf Remark. }\nobreak}{\par \medskip}
\usepackage{fancyhdr}
 \fancypagestyle{plain}{
  \fancyhead{}
  
}

\def\paragraph#1{{\bf #1\ }}
\newcommand{\vertiii}[1]{{\left\vert\kern-0.25ex\left\vert\kern-0.25ex\left\vert #1
    \right\vert\kern-0.25ex\right\vert\kern-0.25ex\right\vert}}
\newcommand{\expo}{\mathrm{e}}

\newcommand{\Var}{\mathrm{Var}}

\newcommand{\dd}{\mathrm{d}}

\newcommand{\HH}{\mathrm{H}}
\newcommand{\overbar}[1]{\mkern 1.5mu\overline{\mkern-1.5mu#1\mkern-1.5mu}\mkern 1.5mu}
\usepackage[utf8]{inputenc}

\def\Proof{\noindent{\bf Proof}\quad}
\def\qed{\hfill$\square$\smallskip}

\title{$K$-averaging agent-based model: propagation of chaos and convergence to equilibrium}

\author{Fei Cao\footnotemark[1]}

\begin{document}
\maketitle

\footnotetext[1]{Arizona State University - School of Mathematical and Statistical Sciences, 900 S Palm Walk, Tempe, AZ 85287-1804, USA}

\tableofcontents

\begin{abstract}
  The paper treats an agent-based model with averaging dynamics to which we refer as the \emph{$K$-averaging model}. Broadly speaking, our model can be added to the growing list of dynamics exhibiting self-organization such as the well-known Vicsek-type models \cite{aldana_phase_2003,aldana_phase_2007,pimentel_intrinsic_2008}. In the $K$-averaging model, each of the $N$ particles updates their position by averaging over $K$ randomly selected particles with additional noise. To make the $K$-averaging dynamics more tractable, we first establish a propagation of chaos type result  in the limit of infinite particle number (i.e. $N \to \infty$) using a martingale technique. Then, we prove the convergence of the limit equation toward a suitable Gaussian distribution in the sense of Wasserstein distance as well as relative entropy. We provide additional numerical simulations to illustrate both results.
\end{abstract}


\section{Introduction}

The collective behavior of various particle systems is a subject of intensive research that has potential applications in biology, physics, economics, and engineering \cite{naldi_mathematical_2010,belmonte_self-propelled_2008,chuang_state_2007}. Different models are proposed to study the emergence of flocking of birds, formation of consensus in opinion dynamics, and phase transitions in network models \cite{motsch_heterophilious_2014,porfiri_effective_2016,chate_modeling_2008,barbaro_phase_2012}. Broadly speaking, all of the aforementioned models are instances of \emph{interacting particle systems}, under various interaction rules among the particles. We refer the readers to \cite{liggett_interacting_2012} for a general introduction into this branch of applied mathematics.

In this work, we investigate a simple model to describe the collective alignment of a group of particles. The model we examine here can be classified in general as an \emph{averaging dynamics} and will be referred to as the \emph{$K$-averaging model}. The readers are encouraged to consult \cite{carlen_kinetic_2013,bertin_boltzmann_2006,bertin_hydrodynamic_2009,boissard_trail_2013} for a variety of models in biology and physics in which averaging plays a key role in the model definition. One important inspiration for the present work is a paper of Maurizio Porfiri and Gil Ariel \cite{porfiri_effective_2016}, which can be thought as a $K$-averaging model on the unit circle. In the $K$-averaging model considered in this manuscript, at each time step, we update the position of each particle (viewed as an element of $\mathbb{R}^d$) according to the average position of its $K$ randomly chosen neighbors while being simultaneously subjected to additive noise (see equation \eqref{dynamics}). Thus, we give the following definition.

\begin{figure}[tbhp]
  \centering
  \includegraphics[width=.9\textwidth]{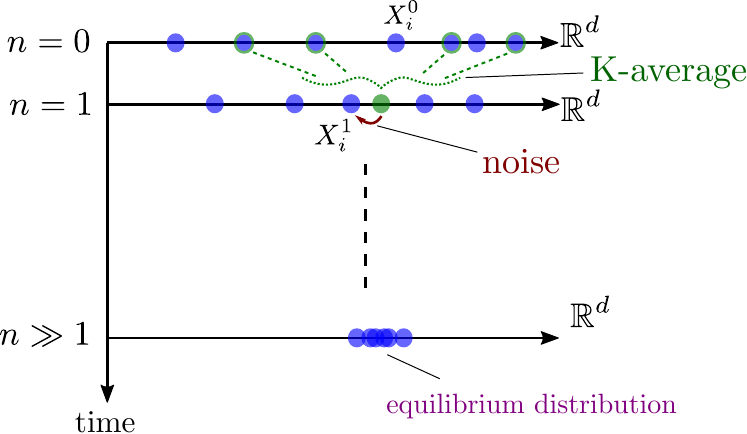}
  \caption{Sketch illustration of the $K$-averaging dynamics \eqref{dynamics}. At each time step, a particle updates its position taking the averaging of $K$ randomly selected particles and add some (Gaussian) noise.}
  \label{fig:illustration}
\end{figure}

\begin{definition}(\textbf{$K$-averaging model})
\label{defK}
Consider a collection of stochastic processes $\{X^n_i\}_{1\leq i\leq N}$ evolving on $\mathbb{R}^d$, where $n$ is the index for time.  At each time step, each particle updates its value to the average of $K$ randomly selected neighbors, subject to an independent noise term:
\begin{equation}\label{dynamics}
X^{n+1}_i := \frac{1}{K}\sum\limits_{j=1}^K X^n_{S^n_i(j)} + W^n_i,\quad 1\leq i\leq N,
\end{equation}
where $S^n_i(j)$ are indices taken randomly from the set $\{1,2,\ldots,N\}$ (i.e., $S^n_i(j) \sim \mathrm{Uniform}(\{1,2,\ldots,N\})$ and is independent of $i,j$ and $n$), and $W^n_i \sim \mathcal{N}(\bm{0},\sigma^2\mathbbm{1}_d)$ is independent of $i$ and $n$, in which $\bm{0}$ and $\mathbbm{1}_d$ stands for the zero vector and the identity matrix in dimension $d$, respectively (see Figure \ref{fig:illustration} for a illustration).
\end{definition}

\begin{figure}[tbhp]
  \centering
  \includegraphics[width=.8\textwidth]{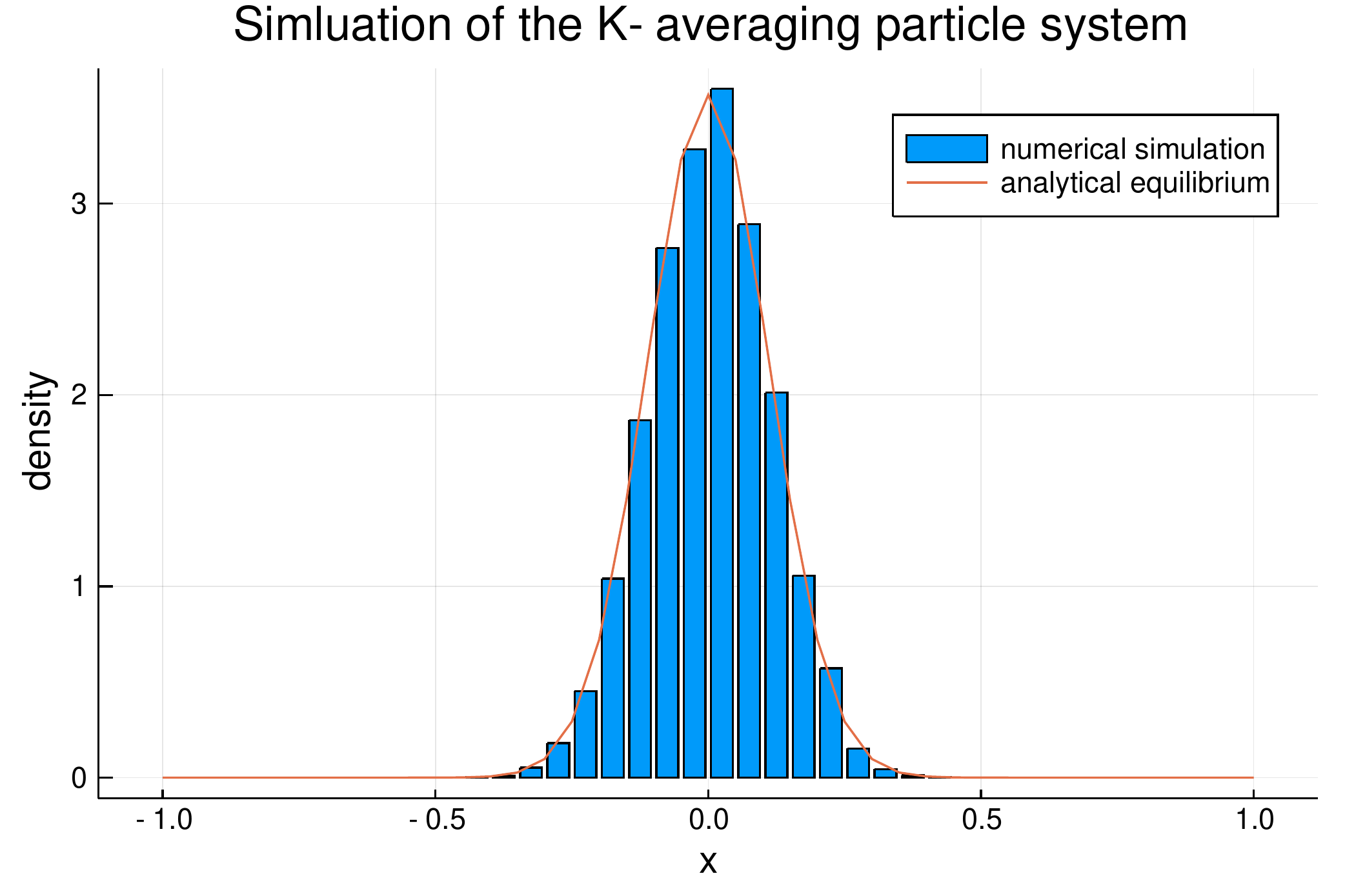}
  \caption{Simulation of the $K$-averaging dynamics in dimension $d = 1$ with $K = 5$ and $N = 5000$ particles after 1000 time steps, in which we used $\sigma = 0.1$ and initially each $X_i \sim \mathrm{Uniform}(-1,1)$. As to be shown later, the distribution of particles will be asymptotically Gaussian under the large $N$ and large time limits.}
  \label{fig:illustration_numerics}
\end{figure}

\begin{figure}[tbhp]
  \centering
  \includegraphics[width=.8\textwidth]{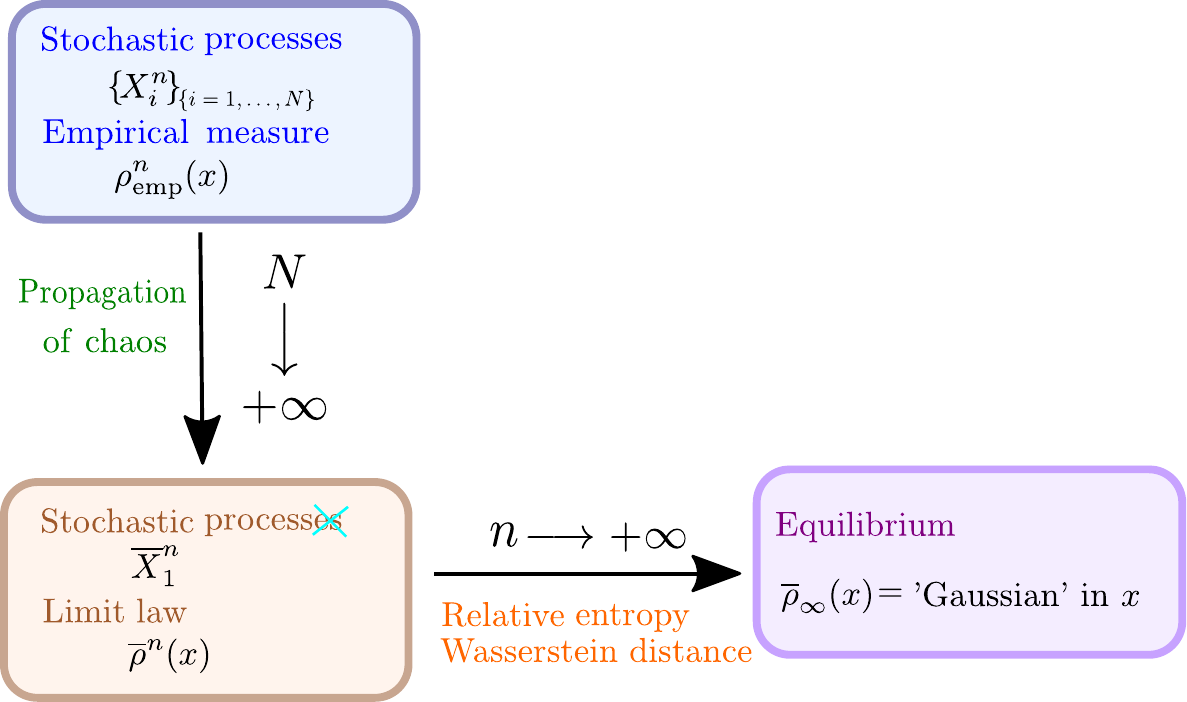}
  \caption{Schematic illustration of the limiting procedure carried out for the study the $K$-averaging dynamics \eqref{dynamics}. The empirical measure $\rho^n_{\mathrm{emp}}(\bm{x})$ of the system (see equation \eqref{emp}) will be shown to converge as $N \to \infty$ to its limit law $\overbar{\rho}^n$ described by the evolution equation \eqref{K-mean-limit}, and then the relaxation of $\overbar{\rho}^n$ to its Gaussian equilibrium will be established.}
  \label{fig:scheme}
\end{figure}

We illustrate the dynamics in Figure \ref{fig:illustration}. The key question of interest is the exploration of the limiting particle distribution as the total number of particles and the number of time steps become large. We illustrate numerically (see Figure \ref{fig:illustration_numerics}) the evolution of the dynamics in dimension $d = 1$ using $N = 5,000$ particles after $n=1000$ time steps.
One of the main difficulty in the rigorous mathematical treatment of models involving large number of interacting particles or agents lies in the general fact interaction will build up correlation over time. Fortunately the framework of kinetic theories allows possible simplification of the analysis of certain such models via suitable asymptotic analysis, see for instance \cite{oelschlager_martingale_1984,meleard_propagation_1987,hauray_n-particles_2007,jabin_quantitative_2018,merle_cutoff_2019,sznitman_topics_1991}.
For the model at hand, our main contribution is two fold: we first prove a result of propagation-of-chaos type under the large $N$ limit (see Theorem \ref{propagationofchaos} for a precise statement), in which interactions among particles are eliminated in finite time and a mean-field dynamics emerges. After the large population limit is carried out and the simplified dynamics (see equation \eqref{limitdyna}) is obtained, we then show that the law of the limiting dynamics defined by \eqref{limitdyna} is asymptotically Gaussian under the large time limit, and such convergence of distribution occurs both in the Waasserstein distance (see Theorem \ref{relaxation}) and in the sense of relative entropy (see Theorem \ref{entropy_relax}). A schematic illustration of the strategy used in this manuscript is presented in Figure \ref{fig:scheme}.
We briefly explain the possible motivation of studying such a model (at least in dimension $d=1$). In the context of a opinion dynamics model (see for instance \cite{baumann_modeling_2020}), $X^n_i$ may represent an evolving opinion of agent $i$ at time step $n$. For a given event, agent $i$ has a opinion $X_i$ (which can be positive or negative) with strength $|X_i|$, and agents update their opinions based on the equation \eqref{dynamics}.

There remain some open questions related to our current work. First, our analysis of the model is restricted to $K\geq 2$, under which we are able to identify the equilibrium distribution and prove various results, we speculate that the propagation of chaos property will be lost if $K=1$, yet we have not been able to find a perfect analytical justification. We also remark here that the case of $K = 1$ can be seen as a variant of the "Choose the Leader" (CL) dynamics introduced in \cite{carlen_kinetic_2013}, in which each of the $N$ particles decides to jump to the location of the other particle chosen independently and uniformly at random at every time step, though noise is injected in such a jump. Second, we think similar results can be obtained if the noise is no longer Gaussian, except that the equilibrium will not be explicit in general.

The remainder of the paper is organized as follows: In section \ref{subsec:2.1}, we present several preliminaries related to random probability measures and the concept of propagation of chaos. Sections \ref{subsec:2.2} and \ref{subsec:2.3} are concerned with the intuitive derivation of the simplified model \eqref{limitdyna} and related properties. We give a full proof of the propagation of chaos result in section \ref{sec:3} and the large time asymptotic of \eqref{limitdyna} is investigated in \ref{sec:4}. We devote section \ref{sec:5} to the continuous-time counterpart of the $K$-averaging model studied in previous sections, and finish the paper with a conclusion in section \ref{sec:6}.

\section{$K$-averaging model}
\label{sec:2}
In section 2.1, we perform a brief review on convergence of random probability measures and the notion of propagation of chaos. Section 2.2 encapsulated a heuristic argument for the large $N$ limit, and we prove a Lipschitz continuity property of the key operator $T$ arising naturally from the $K$-averaging dynamics in section 2.3, which will be leveraged in the proof of Theorem \ref{propagationofchaos}.

\subsection{Review propagation of chaos and convergence of random measures}
\label{subsec:2.1}
We devote this section to a quick review on propagation of chaos and convergence of random probability measures. First, we intend to briefly discuss about propagation of chaos, but we need to carefully define what propagation of chaos means. With this aim, we consider a (stochastic) $N$-particle system denoted by $(X_1,\ldots,X_N)$ in which particles are indistinguishable. In other words, the particle system enjoys a property known as permutation invariance, i.e. for any test function $\varphi$ and permutation $\eta \in \mathcal{S}_N$:
\[\mathbb{E}[\varphi\big(X_1,\ldots,X_N\big)] = \mathbb{E}[\varphi\big(X_{\eta(1)},\ldots,X_{\eta(N)}\big)].\]
In particular, all the single processes $X_i$ have the same law for $1\leq i\leq N$ (although they are in general correlated). Let $\rho^{(N)}(x_1,\ldots,x_N)$ to be the density distribution of the $N$-particle process and denote $\rho^{(N)}_k$ its $k$-particle marginal density, i.e., the law of the process $(X_1,\ldots,X_k)$:

\[\rho^{(N)}_k(x_1,\ldots,x_k) := \int_{x_{k+1},\ldots,x_N} \rho^{(N)}(x_1,\ldots,x_N) \dd x_{k+1}\ldots\dd x_N.\]

Consider now a (potential) limit stochastic process $(\overbar{X}_1,\ldots,\overbar{X}_k)$ where $\{\overbar{X}_i\}_{1\leq i\leq k}$ are i.i.d. Denote by $\overbar{\rho}_1$ the law of a single process, thus by independence assumption the law of all the process is given by:

\[\overbar{\rho}_k(x_1,\ldots,x_k) = \prod\limits_{i=1}^k \overbar{\rho}_1(x_i),\quad \text{i.e.},~\overbar{\rho}_k = \bigotimes_{i=1}^k \overbar{\rho}_1.\]

The following definition is classical and can be found for instance in \cite{sznitman_topics_1991,carlen_kinetic_2013}.
\begin{definition}
We say that the stochastic process $(X_1,\ldots,X_N)$ satisfies the propagation of chaos if for any fixed $k$:
\begin{equation}\label{defn}
\rho^{(N)}_k \xrightharpoonup[]{N \to \infty} \overbar{\rho}_k,
\end{equation}
which is equivalent to the validity of the following relation for any test function $\varphi$:
\begin{equation}\label{equidefn}
\mathbb{E}[\varphi\big(X_1,\ldots,X_k\big)] \xrightarrow[]{N \to \infty} \mathbb{E}[\varphi\big(X_1,\ldots,X_k\big)].
\end{equation}
\end{definition}

Next, we shift to a review on convergence of random probability measures. Such topic can be found for instance in a classical book by Billingsley \cite{billingsley_convergence_2013}. However, we prefer to give a more practical treatment on convergence of random probability measures, based on \cite{berti_almost_2006}. Consider a sequence of random probability measures $\mu_n(\omega)$, i.e., for a given $\omega \in \Omega$, $\mu_n(\omega) \in \mathcal{P}(\mathbb{R}^d)$.  We shall define the mode of convergence as follows:

\begin{definition}
We say that $\mu_n$ converges to $\mu \in \mathcal{P}(\mathbb R^d)$ in probability, denoted by $\mu_n \xrightarrow[]{\mathbb{P}} \mu$, if
\begin{equation}\label{prob}
\langle \mu_n(\omega),\varphi \rangle \xrightarrow[]{\mathbb{P}} \langle \mu(\omega),\varphi \rangle \quad \text{for any $\varphi \in C_b(\mathbb{R}^d)$}.
\end{equation}
\end{definition}

We record here a simple criteria to test the convergence in probability of random measures.

\begin{lemma}
Suppose that the sequence of random measures $\{\mu_n(\omega)\}_n$ satisfies
\begin{equation}\label{useful}
\mathbb{E}_\omega[|\langle \mu_n(\omega)-\mu(\omega),\varphi \rangle|] \xrightarrow[]{n \to \infty} 0 \quad \text{for all $\varphi \in C_b(\mathbb{R}^d)$}.
\end{equation}
Then $\mu_n \xrightarrow[]{\mathbb{P}} \mu$.
\end{lemma}

\Proof It is a direct application of the Markov's inequality. Fixing $\varphi \in C_b(\mathbb{R}^d)$ and let $\varepsilon >0$, we have
\begin{align*}
\mathbb{P}[|\langle \mu_n(\omega),\varphi \rangle - \langle \mu(\omega),\varphi \rangle| > \varepsilon] &= \mathbb{P}[|\langle \mu_n(\omega)-\mu(\omega),\varphi \rangle| > \varepsilon] \\
&\leq \frac{\mathbb{E}_\omega[|\langle \mu_n(\omega)-\mu(\omega),\varphi \rangle|]}{\varepsilon} \xrightarrow[]{n \to \infty} 0.
\end{align*}
Therefore, the random variables $X_n(\omega):=\langle \mu_n(\omega),\varphi \rangle$ converges in probability to $X(\omega):= \langle \mu(\omega),\varphi \rangle$. Since it is true for any $\varphi \in C_b(\mathbb{R}^d)$, we deduce that $\mu_n \xrightarrow[]{\mathbb{P}} \mu$. \qed

\subsection{Formal limit as $N \to \infty$}
\label{subsec:2.2}
We would like to investigate formally the limit as $N \to \infty$ of the dynamics, and we will provide the rigorous derivation in the next section. Motivated by the famous molecular chaos assumption (also known as propagation of chaos), which suggests that we have the statistical independence among the particle systems defined by \eqref{dynamics} under the large $N \to \infty$ limit, we henceforth give the following definition of the limiting dynamics of $\overbar{X}_1$ as $N \to \infty$ from the process point of view.

\begin{definition}(\textbf{Asymptotic $K$-averaging model})
\label{defAK}
We define a collection of random variables $\{\overbar{X}^n\}_{n\geq 0}$ by setting $\overbar{X}^0 = X^0_1$ and
\begin{equation}\label{limitdyna}
\overbar{X}^{n+1} := \frac{1}{K}\sum\limits_{j=1}^K \overbar{Y}^n_{j} + W^n,
\end{equation}
where $\{\overbar{Y}^n_{j}\}_{1\leq j\leq K}$ are $K$ i.i.d copies of $\overbar{X}^n$ and $W^n \sim \mathcal{N}(\bm{0},\sigma^2\mathbbm{1}_d)$ is independent of $n$.
\end{definition}

If we denote $\overbar{\rho}$ to be the law of $\overbar{X}$, then it is possible to determine the evolution of $\overbar{\rho}$ with respect to time $n$. For this purpose, We will first collect some definitions to be used throughout the manuscript.

\begin{definition}\label{def1}
We use $\mathcal{P}(\mathbb{R}^d)$ to represent the space of probability measures on $\mathbb{R}^d$. We will denote by $\phi$ the probability density of a $d$-dimensional Gaussian random variable $\mathcal{E} \sim \mathcal{N}(\bm{0},\sigma^2\mathbbm{1}_d)$. For $\rho \in \mathcal{P}(\mathbb R^d)$, we define $T\colon \mathcal{P}(\mathbb R^d) \to \mathcal{P}(\mathbb R^d)$ through
\begin{equation}\label{operaT}
T[\rho]= \phi*S_K[C_K[\rho]],
\end{equation} in which the $C_K$ is the $K$-fold repeated self-convolution defined via
\begin{equation}\label{operaC}
C_K[\rho] := \underbrace{\rho*\rho*\cdots*\rho}_{\text{$K$ times}},
\end{equation} and $S_K$ is the scaling (renormalization) operator given by
\begin{equation}\label{operaS}
S_K[\rho](\bm{x}):= K^d\cdot\rho(K\bm{x}),\quad \forall \bm{x}\in \mathbb R^d.
\end{equation}
\end{definition}

\begin{remark}
We emphasize here that the operator $T$ given in Definition \ref{def1} fully encodes the update rule \eqref{limitdyna} for the asymptotic $K$-averaging model. Indeed, for each valid test function $\varphi$, we have
\begin{equation}\label{testfunc}
\langle \overbar{\rho}^{n+1},\varphi \rangle = \mathbb{E}[\varphi(\overbar{X}^{n+1})] = \mathbb{E}\left[\varphi\left(\frac{1}{K}\sum\limits_{j=1}^K \overbar{Y}^n_{j} + W^n\right)\right] = \langle T[\overbar{\rho}^n],\varphi \rangle,
\end{equation}
where the last equality follows because the random variable $\frac{1}{K}\sum\limits_{j=1}^K \overbar{Y}^n_{j} + W^n$ has law $T[\overbar{\rho}^n]$. Thus, from the density point of view, as $\overbar{\rho}^n$ is the law of $\overbar{X}$ at time $n$, then $T[\overbar{\rho}^n]$ represents exactly the law of $\overbar{X}$ at time $n+1$.
\end{remark}

Equipped with Definition \ref{def1}, we can write \textbf{the evolution of the limit equation} as
\begin{equation}\label{K-mean-limit}
\overbar{\rho}^{n+1} = T[\overbar{\rho}^n],\quad n\geq 0.
\end{equation}
Notice that the mean value is preserved by the dynamics \eqref{limitdyna}, we will make a harmless assumption throughout this paper that
\begin{equation}\label{normalization}
\int_{\bm{x}\in \mathbb{R}^d} \bm{x}\overbar{\rho}^n(\bm{x}) \dd \bm{x} = \bm{0} ~~\forall n \geq 0.
\end{equation}

\begin{remark}
In dimension 1, we derive from \eqref{limitdyna} that
\[\Var(\overbar{X}^{n+1}) = \Var\left(\frac{1}{K}\sum\limits_{j=1}^K \overbar{Y}^n_{j} + W^n\right) = \frac{\Var(\overbar{X}^n)}{K} + \sigma^2,\]
leading to $\Var(\overbar{X}^n) \xrightarrow[]{n \to \infty} \frac{K\sigma^2}{K-1}$. A similar consideration demonstrates that the covariance matrix associated with $\overbar{X}^n \in \mathbb{R}^d$ converges to $\frac{K\sigma^2}{K-1}\cdot \mathbbm{1}_d$.
\end{remark}

Now we can verify that a suitable Gaussian profile is a fixed point of the iteration process described by \eqref{K-mean-limit} as long as $K\geq 2$.
\begin{lemma}
Fixing $K\geq 2$. Let
\begin{equation}\label{equili}
\overbar{\rho}_\infty(\bm{x}) := \tfrac{1}{(2\pi\sigma^2_\infty)^{\frac{d}{2}}}\expo^{-\frac{|\bm{x}|^2}{2\sigma^2_\infty}}
\end{equation}
with $\sigma^2_\infty := \frac{K}{K-1}\sigma^2$, then $\overbar{\rho}_\infty$ is a fixed point of $T$. i.e, $\overbar{\rho}_\infty$ satisfies $\overbar{\rho}_\infty = T[\overbar{\rho}_\infty]$.
\end{lemma}
\Proof It is readily seen that the operator $T$ maps a Gaussian density to another (possibly different) Gaussian density. We investigate the effect of each operator appearing in the definition of $T$ on $\overbar{\rho}_\infty$. Indeed, since $Z_1 + \cdots + Z_K \sim \mathcal{N}(\bm{0},K\sigma^2_\infty\mathbbm{1}_d)$ when $(Z_i)_{1\leq i\leq K}$ are i.i.d. with law $\mathcal{N}(\bm{0},\sigma^2_\infty\mathbbm{1}_d)$, we have \[C_K[\overbar{\rho}_\infty](\bm{x}) = \tfrac{1}{(2\pi K\sigma^2_\infty)^{\frac{d}{2}}}\expo^{-\frac{|\bm{x}|^2}{2K\sigma^2_\infty}}.\]  Next, notice that $\frac{Z}{K} \sim \mathcal{N}(\bm{0},\frac{\sigma^2_\infty}{K}\mathbbm{1}_d)$ if $Z \sim \mathcal{N}(\bm{0},K\sigma^2_\infty\mathbbm{1}_d)$, from which we deduce that \[S_K[C_K[\overbar{\rho}_\infty]](\bm{x}) = \tfrac{1}{(2\pi \sigma^2_\infty \slash K)^{\frac{d}{2}}}\expo^{-\frac{|\bm{x}|^2}{2\sigma^2_\infty \slash K}}.\] Finally, we conclude that \[T[\overbar{\rho}_\infty](\bm{x}) = \phi*S_K[C_K[\overbar{\rho}_\infty]](\bm{x}) = \tfrac{1}{(2\pi (\sigma^2_\infty \slash K + \sigma^2))^{\frac{d}{2}}}\expo^{-\frac{|\bm{x}|^2}{2(\sigma^2_\infty \slash K+\sigma^2)}} = \overbar{\rho}_\infty(\bm{x}),\] which completes the proof. \qed

\begin{figure}[h]
\centering
\includegraphics[width=.7\textwidth]{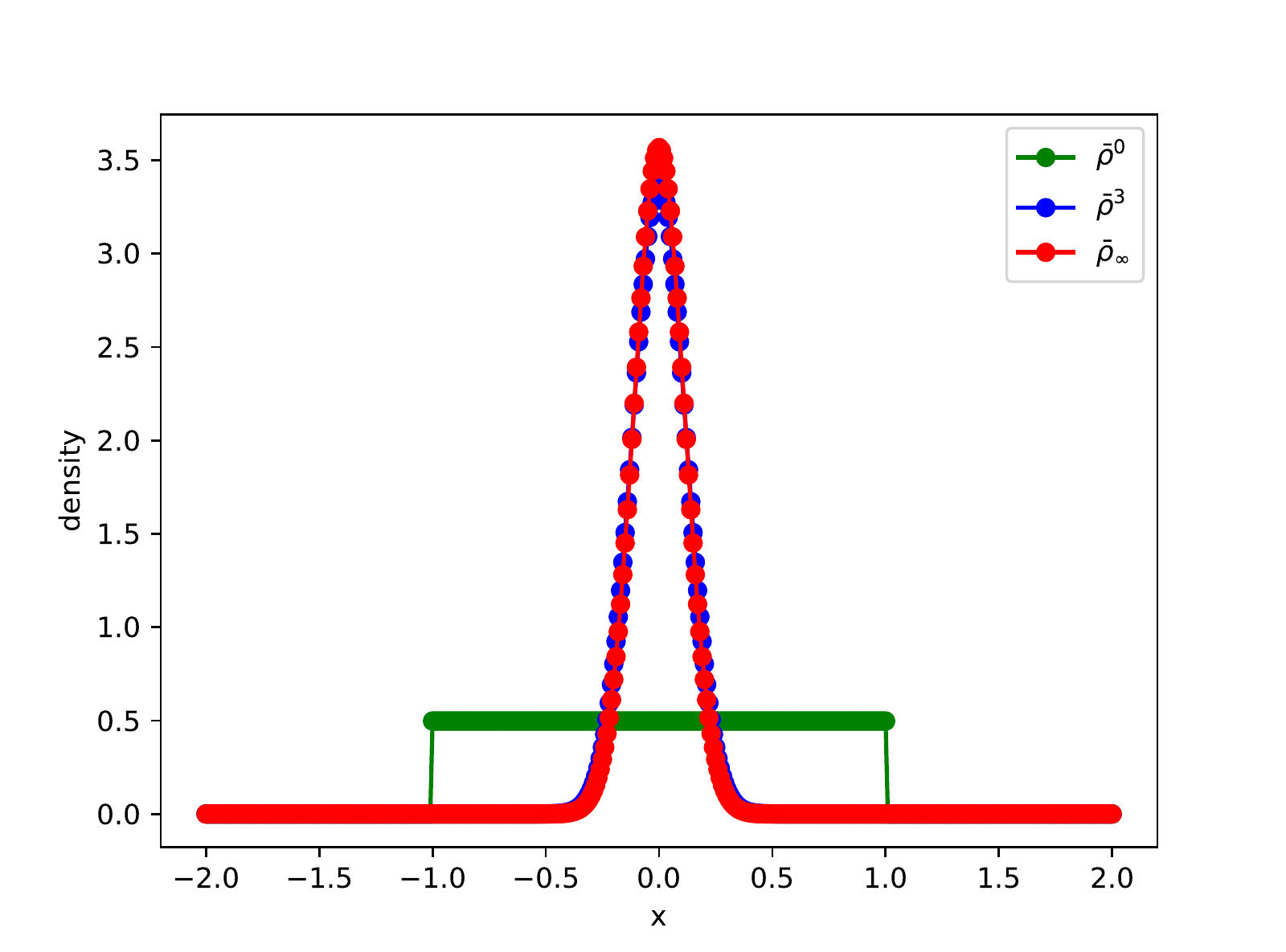}
\caption{Simulation of the discrete evolution equation \eqref{K-mean-limit} in dimension $d = 1$ with $K = 5$ after 3 time steps, in which we used $\sigma = 0.1$ and a uniform distribution over $[-1,1]$ initially $\overbar{\rho}^0(x) := \frac 12\mathbbm{1}_{[-1,1]}(x)$ (the green curve). The blue and red curve represent $\overbar{\rho}^3$ and $\overbar{\rho}_\infty$, respectively. We also remark that in this example $\overbar{\rho}^5$ and $\overbar{\rho}_\infty$ are almost indistinguishable.}
\label{limit_eqn_plot}
\end{figure}

We end this subsection with a numerical experiment demonstrating the relaxation of the solution of \eqref{K-mean-limit} to its Gaussian equilibrium $\overbar{\rho}_\infty$, as is shown in Figure \ref{limit_eqn_plot}.
\subsection{Lipschitz continuity of the operator $T$}
\label{subsec:2.3}
To conclude section 2, we demonstrate a useful property of the operator $T$ introduced in \eqref{operaT}. First, we start with the following definition.

\begin{definition}
For each $\mu \in \mathcal{P}(\mathbb{R}^d)$, we define the (strong) norm of $\mu$, denoted by $\vertiii{\mu}$, via \[\vertiii{\mu} = \sup\limits_{\|\varphi\|_\infty \leq 1} |\langle \mu, \varphi\rangle|.\]
\end{definition}

The main result in this section lies in the Lipschitz continuity of $T$, to which we now turn.

\begin{proposition}\label{LipT}
For each $\mu,\nu \in \mathcal{P}(\mathbb{R}^d)$, we have
\begin{equation}\label{Lipschitz}
\vertiii{T[\mu]-T[\nu]} \leq K\vertiii{\mu-\nu}.
\end{equation}
\end{proposition}

\Proof We recall that for each $g \in \mathcal{P}(\mathbb{R}^d)$ we have \[T[g] = \phi*S_K[C_K[g]].\] Moreover, we have \[\langle S_K[g],h \rangle = K^d\langle g, S_{\frac{1}{K}}[h] \rangle,\quad \forall g,h \in \mathcal{P}(\mathbb{R}^d).\] Also, for $\mu, \nu \in \mathcal{P}(\mathbb{R}^d)$ and $\varphi \in C_b(\mathbb{R}^d)$, there holds
\[\langle \mu*\nu, \varphi\rangle = \langle \nu, \widehat{\mu}*\varphi \rangle,\]
where $\widehat{\mu}$ is defined via $\widehat{\mu}(\bm{x}) := \mu(-\bm{x})$.
Fixing $\varphi$ with $\|\varphi\|_\infty \leq 1$, for each pair of probability measures $\mu,\nu \in \mathcal{P}(\mathbb{R}^d)$, we have
\begin{equation}\label{calc}
\begin{aligned}
\langle T[\mu] - T[\nu], \varphi \rangle &= \langle S_K[C_K[\mu]] - S_K[C_K[\nu]], \phi*\varphi \rangle \\
&= K^d\langle C_K[\mu] - C_K[\nu], S_{\frac{1}{K}}[\phi*\varphi] \rangle \\
&= K^d\sum_{j=0}^{K-1} \langle \underbrace{\mu*\cdots*\mu}_{\text{$j$ times}}*\underbrace{\nu*\cdots*\nu}_{\text{$K-1-j$ times}}*~(\mu-\nu),S_{\frac{1}{K}}[\phi*\varphi] \rangle \\
&:= K^d\sum_{j=0}^{K-1} \langle \kappa_j*(\mu-\nu),S_{\frac{1}{K}}[\phi*\varphi] \rangle \\
&= K^d\langle \mu-\nu, \sum_{j=0}^{K-1} \widehat{\kappa_j}*S_{\frac{1}{K}}[\phi*\varphi]\rangle.
\end{aligned}
\end{equation}
Setting $\psi_j = \widehat{\kappa_j}*S_{\frac{1}{K}}[\phi*\varphi]$ for each $1\leq j\leq K-1$, then we have \[\|\psi_j\|_\infty \leq \|S_{\frac{1}{K}}[\phi*\varphi]\|_\infty \leq \frac{\|\varphi\|_\infty}{K^d} \leq \frac{1}{K^d}.\] Thus, if we define $\varphi^{(1)} = K^d \sum_{j=0}^{K-1} \psi_j$, then $\|\varphi^{(1)}\|_\infty \leq \frac{K^{d+1}}{K^d} = K$. Now taking the supremum over all $\varphi$ with $\|\varphi\|_\infty \leq 1$, we deduce from \eqref{calc} that \[\vertiii{T[\mu]-T[\nu]} \leq K\vertiii{\mu-\nu}\] and the proof is completed. \qed

\section{Propagation of chaos}
\label{sec:3}
This section is devoted to the rigorous proof of propagation of chaos for the $K$-averaging dynamics, by employing a martingale-based technique introduced recently in \cite{merle_cutoff_2019}. We will need the following definition.

\begin{definition}
Let $\{X^n_i\}_{1\leq i\leq N}$ be as in Definition \ref{defK}, we define
\begin{equation}\label{emp}
\rho^n_{\mathrm{emp}}(\bm{x}) := \frac{1}{N}\sum_{i=1}^N \delta_{X^n_i}(\bm{x})
\end{equation} to be the empirical distribution of the system at time $n$. In particular, $\rho^n_{\mathrm{emp}}$ is s stochastic measure.
\end{definition}

Thanks to a classical result (see for instance \textbf{Proposition 1} in \cite{dai_pra_stochastic_2017} or \textbf{Proposition 2.2} in \cite{sznitman_topics_1991}), to justify the propagation of chaos, it suffices to show that \[\rho^n_{\mathrm{emp}} \xrightarrow[]{\mathcal{L}} \overbar{\rho}^n ~\text{as}~ N\to \infty.\] i.e.,
\[\langle \rho^n_{\mathrm{emp}},\varphi \rangle \xrightarrow[]{\mathcal{L}} \langle \overbar{\rho}^n,\varphi \rangle \quad \text{for any $\varphi \in C_b(\mathbb{R}^d)$}.\] In fact, one can prove our first theorem.

\begin{theorem}\label{propagationofchaos}
Under the settings of the $K$-averaging model with $K\geq 2$, if
\begin{equation}\label{ic}
\rho^0_{\mathrm{emp}} \xrightarrow[]{\mathbb{P}} \overbar{\rho}^0 ~\text{as}~ N\to \infty,
\end{equation}
then for each fixed $n \in \mathbb N$ we have \[\rho^n_{\mathrm{emp}} \xrightarrow[]{\mathbb{P}} \overbar{\rho}^n ~\text{as}~ N\to \infty,\] where $\rho^n_{\mathrm{emp}}$ and $\overbar{\rho}^n$ are defined in \eqref{emp} and \eqref{K-mean-limit}, respectively.
\end{theorem}

\Proof We adopt a martingale-based technique developed recently in \cite{merle_cutoff_2019}. We have for each test function $\varphi$ that
\begin{equation}\label{1}
\mathbb{E}\left[\langle \rho^{n+1}_{\mathrm{emp}},\varphi\rangle \right] = \mathbb{E}\left[\frac{1}{N}\sum_{i=1}^N \varphi\left(\frac{1}{K}\sum_{j=1}^K Y^n_{i,j} + W^n_i\right)\right]
\end{equation}

\noindent where $\{Y^n_{i,j}\}$ are i.i.d. with law $\rho^n_{\mathrm{emp}}$. Denoting \[Z^n_i = \frac{1}{K}\sum_{j=1}^K Y^n_{i,j} + W^n_i\] for each $1\leq i\leq N$, since the law of $Z^n_i$ is $T[\rho^n_{\mathrm{emp}}]$ for each $1\leq i\leq n$ and $\{Z^n_i\}_{1\leq i\leq N}$ are i.i.d., following the reasoning behind \eqref{testfunc} we have

\begin{equation}\label{testfunc2}
\mathbb{E}\left[\frac{1}{N}\sum_{i=1}^N \varphi\left(Z^n_i\right)\Bigg|~ \rho^n_{\mathrm{emp}}\right] = \langle T[\rho^n_{\mathrm{emp}}],\varphi\rangle.
\end{equation}

Now if we set
\begin{equation}\label{K-mean}
\begin{aligned}
M^n :&= \frac{1}{N}\sum_{i=1}^N \varphi\left(Z^n_i\right) - \mathbb{E}\left[\frac{1}{N}\sum_{i=1}^N \varphi\left(Z^n_i\right)\Bigg|~ \rho^n_{\mathrm{emp}}\right] \\
&=\langle \rho^{n+1}_{\mathrm{emp}},\varphi \rangle - \langle T[\rho^n_{\mathrm{emp}}],\varphi \rangle
\end{aligned}
\end{equation}
for each $n \geq 0$, then $(M^n)_{n\geq 0}$ defines a martingale. Moreover, thanks to the fact that $\{\varphi(Z^n_i)\}$ are i.i.d. bounded random variables, and using the convention that the variance operation $\Var(X)$ is interpreted as $\Var(X) := \sum_{k=1}^d \Var(X_k)$ when $X$ is a $d$-dimensional vector-valued random variable, we have

\begin{align*}
(\mathbb E[|M^n|])^2 &\leq \mathbb E[|M^n|^2] = \mathbb{E}\left[\Var\left(\frac{1}{N}\sum_{i=1}^N \varphi\left(Z^n_i\right)\Bigg|~ \rho^n_{\mathrm{emp}}\right)\right] \\
&\leq \Var\left(\frac{1}{N}\sum_{i=1}^N \varphi\left(Z^n_i\right)\right) \leq \frac{\|\varphi\|^2_\infty d}{N},
\end{align*}

\noindent where we have employed Popoviciu's inequality (see for instance \cite{popoviciu_sur_1935}) for upper bounding the variance of a bounded random variable. Comparing \eqref{K-mean} with \eqref{K-mean-limit} yields
\begin{equation}\label{keyeqn}
\mathbb{E}\big[|\langle \rho^{n+1}_{\mathrm{emp}} - \overbar{\rho}^{n+1}, \varphi\rangle|\big] \leq \mathbb{E}\big[|\langle T[\rho^n_{\mathrm{emp}}] - T[\overbar{\rho}^n], \varphi\rangle|\big] + \frac{\|\varphi\|^2_\infty d}{\sqrt{N}}.
\end{equation}
Now if $\|\varphi\|_\infty \leq 1$, we can recall the computations carried out in \eqref{calc}, which ensures the existence of some $\varphi^{(1)}$ with $\|\varphi^{(1)}\|_\infty \leq K$ such that
\begin{equation}\label{1ststep}
\langle T[\rho^n_{\mathrm{emp}}] - T[\overbar{\rho}^n], \varphi\rangle = \langle \rho^n_{\mathrm{emp}} - \overbar{\rho}^n, \varphi^{(1)}\rangle.
\end{equation}

\noindent Then we can deduce from \eqref{keyeqn} and \eqref{1ststep} that
\begin{equation}\label{iter}
\mathbb{E}\big[|\langle \rho^{n+1}_{\mathrm{emp}} - \overbar{\rho}^{n+1}, \varphi\rangle|\big] \leq \mathbb{E}\big[|\langle \rho^n_{\mathrm{emp}}) - \overbar{\rho}^n, \varphi^{(1)}\rangle|\big] + \frac{d}{\sqrt{N}},
\end{equation}

\noindent in which $\varphi^{(1)}$ satisfies $\|\varphi^{(1)}\|_\infty \leq K$. We can iterate \eqref{iter} to arrive at
\begin{equation}\label{goal}
\mathbb{E}\big[|\langle \rho^n_{\mathrm{emp}} - \overbar{\rho}^n, \varphi\rangle|\big] \leq \mathbb{E}\big[|\langle \rho^0_{\mathrm{emp}} - \overbar{\rho}^0, \varphi^{(n)} \rangle|\big] + \frac{dn}{\sqrt{N}},
\end{equation}

\noindent in which $\varphi^{(n)}$ satisfies $\|\varphi^{(n)}\|_\infty \leq K^n$. Finally, combining \eqref{ic} with \eqref{goal} allows us to complete the proof of Theorem \ref{propagationofchaos}.\qed

\begin{remark}
As we do not have a uniform-in-time propagation of chaos, we would like to know whether the convergence declared in Theorem \ref{propagationofchaos} still holds if we do not fix $n$ (i.e., if $n \to \infty$). We speculate such an uniform in time convergence can no longer be hoped for by looking at the evolution of the \emph{center of mass} of the particle systems. Indeed, define \[\mathcal{C}_{n+1} := \frac{1}{N}\sum_{i=1}^N X^n_i\] to be the location of the center of mass, and denote by $\mathcal{F}_n$ the natural filtration generated by $(X^n_1,\cdots,X^n_N)$, then in dimension $d = 1$ we have

\begin{align*}
\mathbb{E}[\mathcal{C}_{n+1} \mid \mathcal{F}_n] &= \frac{1}{N}\sum_{i=1}^N \mathbb{E}[X^{n+1}_i \mid \mathcal{F}_n] = \frac{1}{N}\sum_{i=1}^N \frac{1}{K}\sum_{j=1}^K \mathbb{E}[X^n_{S^n_i(j)} \mid \mathcal{F}_n] \\
&= \frac{1}{N}\sum_{i=1}^N \frac{1}{K}\sum_{j=1}^K \frac{1}{N}\sum_{\ell= 1}^N X^n_\ell = \mathcal{C}_n,
\end{align*}

\begin{align*}
\hspace{-0.65in} \mathbb{E}[(\mathcal{C}_{n+1}-\mathcal{C}_n)^2 \mid \mathcal{F}_n] &= \mathbb{E}\left[\left(\frac{1}{N}\sum_{i=1}^N X^{n+1}_i-\mathcal{C}_n\right)^2 \Bigg|~ \mathcal{F}_n\right] \\
                                                                                    &= \Var\left[\frac{1}{N}\sum_{i=1}^N X^{n+1}_i \Bigg|~ \mathcal{F}_n \right] = \frac{1}{N}\Var[X^{n+1}_1 \mid \mathcal{F}_n] \\
  &\geq \frac{\sigma^2}{N},
\end{align*}

\noindent where the last equality comes from the fact that $X^{n+1}_i$ and $X^{n+1}_j$ are i.i.d. \emph{given $\mathcal{F}_n$}. Thus, loosely speaking, at least in dimension $d = 1$, the center of mass of the particle systems behaves like a discrete time Brownian motion with intensity of order at least $\mathcal{O}(1\slash \sqrt{N})$, such an variation can accumulate in time which will eventually ruin the chaos propagation property in the long run.
\end{remark}

\section{Large time behavior}
\label{sec:4}
The long time behavior of the limit equation, resulted from the simplified mean-field dynamics, is treated in this section. In section 4.1, by employing a coupling technique and equipping the space of probability measures on $\mathbb R^d$ with the Wasserstein distance, we will justify the asymptotic Gaussianity of the distribution of each particle. Then we will strengthen the convergence result shown in the previous section in section 4.2, and numerical simulations are also performed in support of our theoretical discoveries in section 4.3. We emphasize here that coupling techniques will be at the core of our proof in section 4.1, and the technique used in section 4.2 depends heavily on several classical results in information theory.

\subsection{Convergence in Wasserstein distance}
\label{subsec:4.1}
After we have achieved the transition from the interacting particle system \eqref{dynamics} to the simplified de-coupled dynamics \eqref{limitdyna} under the limit $N\to \infty$, in this section we will analyze \eqref{limitdyna} and its associated evolution of its law (governed by \eqref{K-mean-limit}), with the intention of proving the convergence of $\overbar{\rho}^n$ to a suitable Gaussian density. The main ingredient underlying our proof lies in a coupling technique. First, we recall the following classical definition.

\begin{definition}
The Wasserstein distance (of order 2) is defined via\[\mathcal{W}^2_2(\mu,\nu):= \inf_{\substack{X\sim \mu \\ Y\sim \nu}} \mathbb{E}[|X-Y|^2],\] where both $\mu$ and $\nu$ are probability measures on $\mathbb{R}^d$.
\end{definition}

\noindent We can now state and prove our main result in this section.
\begin{theorem}\label{relaxation}
Assume that the innocent-looking normalization \eqref{normalization} holds and $K\geq 2$, then for the dynamics \eqref{K-mean-limit}, we have
\begin{equation}\label{contraction}
\mathcal{W}^2_2(\overbar{\rho}^{n+1},\overbar{\rho}_\infty) \leq \frac{1}{K}\mathcal{W}^2_2(\overbar{\rho}^n,\overbar{\rho}_\infty),~\forall n\geq 0
\end{equation}
In particular, if $\overbar{\rho}^0 \in \mathcal{P}(\mathbb R^d)$ is chosen such that $\mathcal{W}^2_2(\overbar{\rho}^0,\overbar{\rho}_\infty) < \infty$, then \[\lim_{n\to\infty} \mathcal{W}^2_2(\overbar{\rho}^n,\overbar{\rho}_\infty) = 0.\]
\end{theorem}

\Proof We first show that
\begin{equation}\label{Wasser}
\mathcal{W}^2_2(T(\mu),T(\nu)) \leq \frac{1}{K}\mathcal{W}^2_2(\mu,\nu)
\end{equation}
for each $\mu,\nu \in \mathcal{P}(\mathbb R^d)$. In other words, if we equip the space $\mathcal{P}(\mathbb R^d)$ with the Wasserstein distance of order 2, $T$ is a strict contraction as long as $K\geq 2$. Now we fix $\mu,\nu \in \mathcal{P}(\mathbb R^d)$. It is recalled that $T(\mu)$ is the law of the random variable \[X := \frac{X_1+\cdots+X_K}{K} + \mathcal{E},\] where $\{X_i\}_{1\leq i\leq K}$ are i.i.d. with law $\mu$ and $\mathcal{E} \sim \mathcal{N}(\bm{0},\sigma^2\mathbbm{1}_d)$. Thus, if we also introduce \[Y := \frac{Y_1+\cdots+Y_K}{K} + \widetilde{\mathcal{E}},\] in which $\{Y_i\}_{1\leq i\leq K}$ are i.i.d. with law $\nu$ and $\widetilde{\mathcal{E}} \sim \mathcal{N}(\bm{0},\sigma^2\mathbbm{1}_d)$, then we can write
\begin{align*}
\mathcal{W}_2^2(T(\mu),T(\nu)) &= \inf_{X\sim T(\mu),Y\sim T(\nu)} \mathbb{E}[|X-Y|^2] \\
&= \inf_{X_i\sim \mu,Y_i\sim \nu} \mathbb{E}\left[\left|\frac{X_1+\cdots+X_K}{K}+\mathcal{E}-\frac{Y_1+\cdots+Y_K}{K}-\widetilde{\mathcal{E}}\right|^2\right].
\end{align*}
We can couple $(X_1,\cdots,X_K,\mathcal{E})$ and $(Y_1,\cdots,Y_k,\tilde{\mathcal{E}})$ as we want. First, we take $\mathcal{E}=\widetilde{\mathcal{E}}$, meaning we have a common source of noise. Second, fixing $\eta>0$, we take $(X_1,Y_1)$ such that \[\mathbb{E}[|X_1-Y_1|^2] \leq \mathcal{W}_2^2(\mu,\nu)+\eta,\] (i.e. almost best coupling). Finally, we perform similarly for the other $(X_i,Y_i)$ with $(X_i,Y_i)$ independent of $(X_j,Y_j)$ if $i\neq j$. These procedures lead us to
\begin{align*}
\mathcal{W}_2^2(T(\mu),T(\nu)) &\leq \mathbb{E}\left[\left|\frac{X_1-Y_1+\cdots+X_K-Y_K}{K}\right|^2\right] \\
&\leq \frac{1}{K^2}\left(\mathbb{E}[|X_1-Y_1|^2] + \cdots + \mathbb{E}[|X_K-Y_K|^2]\right) \\
&\leq \frac{1}{K}\mathcal{W}_2^2(\mu,\nu)+ \frac{\eta}{K}.
\end{align*}
Since this is true for any $\eta>0$, \eqref{Wasser} is verified. Now we can deduce from \eqref{Wasser} that
\[\mathcal{W}^2_2(\overbar{\rho}^{n+1},\overbar{\rho}_\infty) = \mathcal{W}^2_2(T(\overbar{\rho}^n),T(\overbar{\rho}_\infty)) \leq \frac{1}{K}\mathcal{W}^2_2(\overbar{\rho}^n,\overbar{\rho}_\infty),\] whence \eqref{contraction} is proved. \qed

\subsection{Convergence in relative entropy}
\label{subsec:4.2}
In this subsection we will show that the evolution of the discrete equation \eqref{K-mean-limit} relaxes to its Gaussian equilibrium $\overbar{\rho}_\infty$ in the sense of relative entropy, as long as $K\geq 2$. Before stating our result, we first clarify some definitions. We refer the reader to \cite{cover_elements_1999} for a comprehensive account of modern information theory.

\begin{definition}
\label{def_entropy}
We use \[\HH(X):=\HH(g) = \int_{\mathbb R^d} g(\bm{x})\log g(\bm{x}) \dd \bm{x}\] to represent the differential entropy of a $\mathbb R^d$-valued random variable $X$ with law $g$. Moreover, \[\mathrm{D}_{\mathrm{KL}}(g || h) := \HH(g) - \HH(g,h) = \int_{\mathbb R^d} g(\bm{x})\log g(\bm{x}) \dd \bm{x}- \int_{\mathbb R^d} g(\bm{x})\log h(\bm{x}) \dd \bm{x}\] denotes the relative entropy from $h \in \mathcal P(\mathbb R^d)$ to $g \in \mathcal P(\mathbb R^d)$, in which \[\HH(g,h) := \HH(X,Y) = \int_{\mathbb R^d} g(\bm{x})\log h(\bm{x}) \dd \bm{x}\] is the cross-entropy from $g$ to $h$ (or equivalently, from $X$ to $Y$ where the laws of $X$ and $Y$ are $g$ and $h$, respectively).
\end{definition}

For the reader's convenience, we explicitly state two fundamental results from information theory that we shall reply on.

\begin{lemma}[Shannon-Stam]
\label{lem:Shannon-Stam}
Under the set-up of Definition \ref{def_entropy}, we have \[\HH(\sqrt{\lambda}X + \sqrt{1-\lambda}Y) \leq \lambda\HH(X) + (1-\lambda)\HH(Y)\] for each $\lambda \in [0,1]$.
\end{lemma}

Lemma \ref{lem:Shannon-Stam} is one of the three equivalent formulations of the well-known Shannon-Stam inequality, see for instance section 1.3.2 of \cite{rezakhanlou_entropy_2008}. The next lemma (see for instance \textbf{Theorem 1} in \cite{artstein_solution_2004-1} or equation (7) in \cite{madiman_generalized_2007}) demonstrates the monotonicity of the differential entropy along re-scaled sum of i.i.d. square-integrable random variables.

\begin{lemma}
\label{lem:monotonicity}
Let $X_1,X_2,\ldots$ be i.i.d. square-integrable random variables. Then \[\HH\left(\frac{X_1+\cdots+X_n}{\sqrt{n}}\right) \leq \HH\left(\frac{X_1+\cdots+X_{n-1}}{\sqrt{n-1}}\right)\] for each $n \geq 2$.
\end{lemma}

\begin{theorem}\label{entropy_relax}
Assume that $\overbar{\rho}$ is a solution to \eqref{K-mean-limit}, then for each fixed $K \geq 2$ we have
\begin{equation}\label{discrete}
\mathrm{D}_{\mathrm{KL}}(\overbar{\rho}^{n+1} || \overbar{\rho}_\infty) \leq \frac{1}{K}\mathrm{D}_{\mathrm{KL}}(\overbar{\rho}^n || \overbar{\rho}_\infty).
\end{equation}
In particular, for each $K \geq 2$ we have $\mathrm{D}_{\mathrm{KL}}(\overbar{\rho}^n || \overbar{\rho}_\infty) \to 0$ as $n \to \infty$.
\end{theorem}

\Proof Let $\{\overbar{X}^n\}_{n\geq 0}$ be as in Definition \ref{defAK}. If we introduce a random variable $\overbar{X}_\infty$ with law $\overbar{\rho}_\infty$, i.e.,$\overbar{X}_\infty \sim \mathcal{N}(\bm{0},\sigma^2_\infty \mathbbm{1}_d)$, then for each $n \in \mathbb{N}$, we can rewrite \eqref{limitdyna} as
\[\overbar{X}^{n+1} = \frac{1}{\sqrt{K}}\cdot \frac{1}{\sqrt{K}}\sum_{j=1}^K \overbar{Y}^n_j + \sqrt{\frac{K-1}{K}}\overbar{X}_\infty,\]
since $\sqrt{\frac{K-1}{K}}\overbar{X}_\infty = W^n$ in law. Setting $\gamma = \frac{1}{\sqrt{K}}$, we obtain
\[\overbar{X}^{n+1} = \sqrt{\gamma}\cdot \frac{1}{\sqrt{K}}\sum_{j=1}^K \overbar{Y}^n_j + \sqrt{1-\gamma}\cdot\overbar{X}_\infty.\] Consequently, the Shannon-Stam inequality (see Lemma \ref{lem:Shannon-Stam}) together with the monotonicity of differential entropy along normalized sum of i.i.d. random variables (see Lemma \ref{lem:monotonicity}) yields
\begin{equation}\label{iden1}
\HH(\overbar{X}^{n+1}) \leq \gamma\HH\left(\frac{1}{\sqrt{K}}\sum_{j=1}^K \overbar{Y}^n_j\right) + (1-\gamma)\HH(\overbar{X}_\infty) \leq \gamma\HH(\overbar{X}^n) + (1-\gamma)\HH(\overbar{X}_\infty).
\end{equation}
Next, we observe that the cross-entropy from each $f \in \mathcal{P}(\mathbb R^d)$ with mean $\bm{0}$ to the equilibrium distribution $\overbar{\rho}_\infty$ is essentially the variance of $f$, meaning that \[\HH(f,\overbar{\rho}_\infty) = -\frac{d}{2}\log(2\pi\sigma^2_\infty) - \frac{\int_{\mathbb R^d} |\bm{x}|^2f(\bm{x}) \dd \bm{x}}{2\sigma^2_\infty}.\]
In particular, if $X$ and $Y$ are independent random variables with mean $\bm{0}$ and $a^2+b^2 = 1$, then
\[\HH(aX+bY,\overbar{X}_\infty) = a^2\HH(X,\overbar{X}_\infty) + b^2\HH(Y,\overbar{X}_\infty).\]
Thus, using this formulation with $\sqrt{\gamma}$ and $\sqrt{1-\gamma}$, we find
\begin{align*}
\HH(\overbar{X}^{n+1},\overbar{X}_\infty) &= \gamma\HH\left(\frac{1}{\sqrt{K}}\sum_{j=1}^K \overbar{Y}^n_j,\overbar{X}_\infty\right) + (1-\gamma)\HH(\overbar{X}_\infty,\overbar{X}_\infty) \\
&= \gamma\HH(\overbar{X}^n,\overbar{X}_\infty) + (1-\gamma)\HH(\overbar{X}_\infty).
\end{align*}
Combining this with \eqref{iden1} leads to
\begin{align*}
\mathrm{D}_{\mathrm{KL}}(\overbar{\rho}^{n+1} || \overbar{\rho}_\infty) &= \HH(\overbar{X}^{n+1}) - \HH(\overbar{X}^{n+1},\overbar{X}_\infty) \\
&\leq \gamma\HH(\overbar{X}^n) + (1-\gamma)\HH(\overbar{X}_\infty) - \gamma\HH(\overbar{X}^n,\overbar{X}_\infty) - (1-\gamma)\HH(\overbar{X}_\infty) \\
&= \gamma\mathrm{D}_{\mathrm{KL}}(\overbar{\rho}^n || \overbar{\rho}_\infty),
\end{align*}
and the proof is completed. \qed

\begin{remark}
By Talagrand's inequality (see for instance \textbf{Theorem 9.2.1} in \cite{bakry_analysis_2013}), the convergence $\mathrm{D}_{\mathrm{KL}}(\overbar{\rho}^n || \overbar{\rho}_\infty) \to 0$ implies the convergence $\mathcal{W}^2_2(\overbar{\rho}^n,\overbar{\rho}_\infty) \to 0$.
\end{remark}

\subsection{Numerical illustration of decay in relatively entropy}
\label{subsec:4.3}
\begin{figure}[!htb]
\centering
\includegraphics[scale = 0.6]{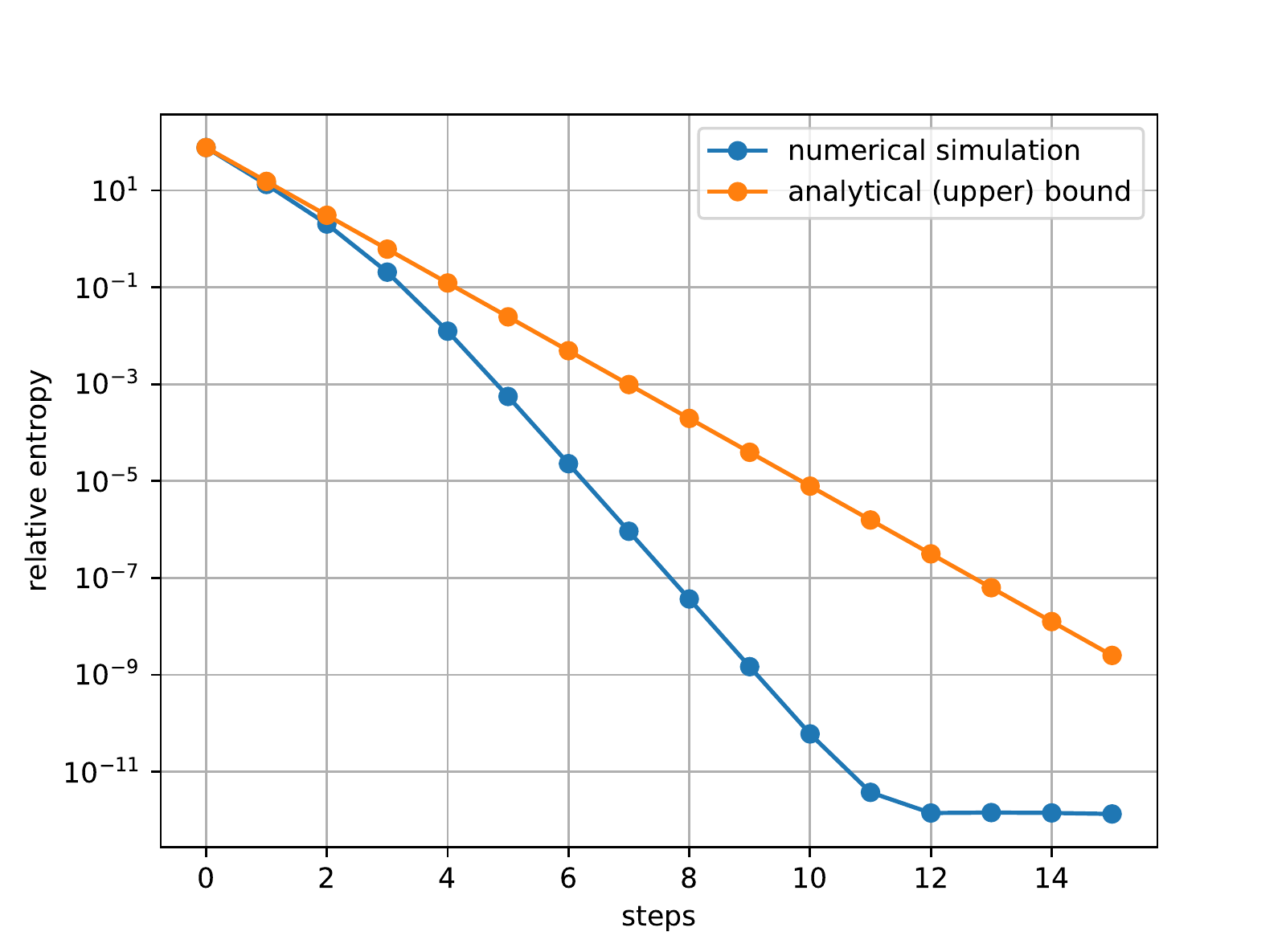}
\caption{Simulation of the relative entropy from $\overbar{\rho}$ to $\overbar{\rho}_\infty$ in dimension $d = 1$ with $K = 5$ after 15 time steps, in which we used $\sigma = 0.1$ and a Laplace distribution $\overbar{\rho}^0(x) = \frac 12\expo^{-|x|}$ initially. The blue and orange curve represent the numerical error and the analytical upper bound on the error, respectively. We also noticed that the numerical error can not really go below $10^{-12}$, but this is presumably due to the floating-point precision error.}
\label{entropy_decay}
\end{figure}

We investigate numerically the convergence of the solution $\overbar{\rho}^n$ of \eqref{K-mean-limit} to its equilibrium $\overbar{\rho}$ in support of our Theorem \ref{entropy_relax}, see Figure \ref{entropy_decay}. We use $d = 1$ (dimension), $K = 5$ (number of neighbors to be averaged over), $\sigma = 0.1$ (the intensity of a centered Gaussian noise) in the simulation of the evolution equation \eqref{K-mean-limit}. To discretize \eqref{K-mean-limit}, we employ the step-size $\Delta x = 0.001$ and a cutoff threshold $M = 100,000$ so that the support of $\overbar{\rho}^n$ is contained in $\{j\Delta x\}_{-M\leq j\leq M}$ for all $n$, and the total number of simulation steps is set to 15. As initial condition, we use the Laplace distribution $\overbar{\rho}^0(x) = \frac 12\expo^{-|x|}$. Moreover, the simulation result is displayed in the semi-logarithmic scale, which clearly indicates a geometrically fast convergence.

\section{Continuous-time $K$-averaging dynamics}
\label{sec:5}
With suitable modifications the argument used in the discrete-time applies in continuous-time as well, so in this section we briefly consider the continuous version of the $K$-averaging model studied in previous sections, i.e., the $K$-averaging occurs according to a Poisson process. First, we give a formal definition of the model.

\begin{definition}(\textbf{Continuous-time $K$-averaging model})
Consider a collection of stochastic processes $\{X_i(t)\}_{1\leq i\leq N}$ evolving on $\mathbb{R}^d$. At each time a Poisson clock with rate $\lambda$ rings, we pick a particle $i \in \{1,\ldots,N\}$ uniformly at random and update the position of $X_i$ according to the average position of $K$ randomly selected neighbors, subject to an independent noise term, i.e., for each test function $\varphi$ the process must satisfy
\begin{equation}\label{contdynamics}
\begin{aligned}
\dd \mathbb{E}[\varphi\big(X_1(t),\ldots,X_N(t)\big)] = \lambda\sum_{i=1}^N \mathbb{E}[&\varphi\big(X_1(t),\ldots,Z_i(t),\ldots,X_N(t)\big) \\
& -\varphi\big(X_1(t),\ldots,X_N(t)\big)]\dd t,
\end{aligned}
\end{equation}
where $Z_i(t) := \frac{1}{K}\sum_{j=1}^K X_{S_i(j)}(t) + W_i(t)$, $S_i(j)$ are indices taken randomly from the set $\{1,2,\ldots,N\}$ (i.e., $S_i(j) \sim \mathrm{Uniform}(\{1,2,\ldots,N\})$ and is independent of $i,j$ and $t$), and $W_i(t) \sim \mathcal{N}(\bm{0},\sigma^2\mathbbm{1}_d)$ is independent of $i$ and $t$.
\end{definition}

In the large $N$ limit, we expect an emergence of a simplified dynamics, which motivates the following definition.

\begin{definition}(\textbf{Asymptotic continuous-time $K$-averaging model})
Consider a $\mathbb{R}^d$-valued stochastic process $\overbar{X}(t)$ which satisfies the following relation for each test function $\varphi$:
\begin{equation}\label{ContK}
\dd \mathbb{E}[\varphi\big(\overbar{X}(t)\big)] = \lambda\mathbb{E}[\varphi\big(\overbar{Z}(t)\big) - \varphi\big(\overbar{X}(t)\big)]\dd t,
\end{equation}
in which $\overbar{Z}(t) := \frac{1}{K}\sum_{j=1}^K \overbar{Y}_j(t) + W(t)$, where $\{\overbar{Y}_j(t)\}_{1\leq j\leq K}$ are $K$ i.i.d. copies of $\overbar{X}(t)$ and $W(t) \sim \mathcal{N}(\bm{0},\sigma^2\mathbbm{1}_d)$ is independent of $t$.
\end{definition}

If we define $\overbar{\rho}(\bm{x},t)$ to be the law of $\overbar{X}$ at time $t$, one can readily see that the evolution of $\overbar{\rho}$ is governed by

\begin{equation}\label{contPDE}
\partial_t \overbar{\rho} = \lambda(T[\overbar{\rho}] - \overbar{\rho}),~~t\geq 0.
\end{equation}

Moreover, one can show the continuous-time analog of Theorem \ref{propagationofchaos} and Theorem \ref{entropy_relax}.

\begin{theorem}\label{continuous_counterpart}
Let $\rho_{\mathrm{emp}}(t) := \frac{1}{N}\sum_{i=1}^N \delta_{X_i(t)}$ to be the empirical distribution of the
system determined by \eqref{contdynamics} at time $t$ and $\overbar{\rho}$ the solution of \eqref{contPDE} with the Gaussian equilibrium $\overbar{\rho}_\infty$ defined in \eqref{equili}, then
\begin{enumerate}[label=(\roman*)]
\item under the set-up of the continuous-time $K$-averaging model with $K \geq 2$, if
\begin{equation}\label{assump}
\rho_{\mathrm{emp}}(0) \xrightarrow[]{\mathbb{P}} \overbar{\rho}(0) ~\text{as}~ N\to \infty,
\end{equation}

then we have \[\rho_{\mathrm{emp}}(t) \xrightarrow[]{\mathbb{P}} \overbar{\rho}(t) ~\text{as}~ N\to \infty,\] holding for all $0\leq t \leq T$ with any prefixed $T >0$.
\item for each fixed $K \geq 2$ we have
\begin{equation}\label{entdec}
\frac{\dd}{\dd t}\mathrm{D}_{\mathrm{KL}}(\overbar{\rho}||\overbar{\rho}_\infty) \leq -\lambda(1-\gamma) \mathrm{D}_{\mathrm{KL}}(\overbar{\rho}||\overbar{\rho}_\infty),
\end{equation}
where $\gamma = \frac{1}{K}$ as before. In particular, we have
\begin{equation}\label{conseq}
\mathrm{D}_{\mathrm{KL}}(\overbar{\rho}(t) || \overbar{\rho}_\infty) \leq \mathrm{D}_{\mathrm{KL}}(\overbar{\rho}(0) || \overbar{\rho}_\infty)\cdot \expo^{-\lambda(1-\gamma)t}.
\end{equation}
\end{enumerate}
\end{theorem}

\Proof We assume without loss of generality that $\lambda =1$. For $(i)$, mimic the argument in the discrete-time setting we obtain for each test function $\varphi$ that
\begin{align*}
\dd \mathbb{E}\left[\langle \rho_{\mathrm{emp}}(t), \varphi \rangle \right] &= \mathbb{E}\left[\frac{1}{N}\sum_{i=1}^N \varphi\big(Z_i(t)\big) - \frac{1}{N}\sum_{i=1}^N \varphi\big(X_i(t)\big) \right] \dd t \\
&= \mathbb{E}\left[\langle T[\rho_{\mathrm{emp}}(t)] - \rho_{\mathrm{emp}}(t), \varphi \rangle \right] \dd t,
\end{align*}
in which $Z_i(t) = \frac{1}{K}\sum_{j=1}^K Y_{i,j}(t) + W_i(t)$ and $\{Y_{i,j}(t)\}$ are i.i.d. with law $\rho_{\mathrm{emp}}(t)$. Then by Dynkin's formula, the compensated process \[M_\varphi(t) := \langle \rho_{\mathrm{emp}}(t), \varphi \rangle - \langle \rho_{\mathrm{emp}}(0), \varphi \rangle - \int_0^t \langle T[\rho_{\mathrm{emp}}(s)] - \rho_{\mathrm{emp}}(s), \varphi \rangle \dd s \] defines a martingale. Comparing with \eqref{contPDE} yields
\begin{equation}\label{error}
\begin{aligned}
\langle \rho_{\mathrm{emp}}(t)-\overbar{\rho}(t), \varphi \rangle| &\leq |M_\varphi(t)| + |\langle \rho_{\mathrm{emp}}(0)-\overbar{\rho}(0), \varphi \rangle| \\
&+ \int_0^t |\langle T[\overbar{\rho}(s)] - T[\rho_{\mathrm{emp}}(s)] - \big(\overbar{\rho}(s)-\rho_{\mathrm{emp}}(s)\big), \varphi \rangle| \dd s.
\end{aligned}
\end{equation}

We then take the supremum over all $\varphi$ with $\|\varphi\|_\infty \leq 1$ to deduce from Proposition \ref{LipT} and \eqref{error} that
\begin{equation*}
\vertiii{\rho_{\mathrm{emp}}(t)-\overbar{\rho}(t)} \leq \eta(t) + (K+1)\int_0^t \vertiii{\rho_{\mathrm{emp}}(s)-\overbar{\rho}(s)} \dd s,
\end{equation*}

\noindent where we have set \[\eta(t) := \sup\limits_{\|\varphi\|_\infty \leq 1} |M_\varphi(t)| + \vertiii{\rho_{\mathrm{emp}}(0)-\overbar{\rho}(0)}.\] By Gronwall's inequality, we obtain
\begin{equation*}
\sup\limits_{t\in [0,T]} \vertiii{\rho_{\mathrm{emp}}(t)-\overbar{\rho}(t)} \leq \left(\sup\limits_{t\in [0,T]} \eta(t)\right)\expo^{(K+1)T}.
\end{equation*}

\noindent In order to justify our claim $(i)$ for $t\leq T$, it therefore suffices to show that
\begin{equation}\label{err}
\sup\limits_{t\in [0,T]} \eta(t) \xrightarrow[N \to \infty]{\mathbb{P}} 0.
\end{equation}

To show \eqref{err}, we address each term appearing in the definition of $\eta(t)$ separately. The second one vanishes due to our assumption \eqref{assump}. For the first one, i.e., the martingale term, we note that the $i$-th coordinate of $M_\varphi$ is a continuous time martingale with jumps of size $\frac{1}{N}\varphi(Z_i) - \frac{1}{N}\varphi(X_i)$ whose rates of occurrence are $\lambda\cdot \dd t = \dd t$. Therefore,
\begin{equation*}
\mathbb{E}[|M_\varphi(T)|^2] \leq \int_0^T \mathbb{E}\left[\sum_{i=1}^N \left|\frac{1}{N}\varphi(Z_i) - \frac{1}{N}\varphi(X_i)\right|^2\right]\dd t \leq \frac{4\|\varphi\|^2_\infty}{N}T \leq \frac{4T}{N},
\end{equation*}
whence the convergence \[\sup\limits_{t\in [0,T]}\left(\sup\limits_{\|\varphi\|_\infty \leq 1} |M_\varphi(t)|\right) \xrightarrow[N \to \infty]{\mathbb{P}} 0 \] follows readily from Doob's martingale inequality. For $(ii)$, we recall that in the discrete-time case (with $\gamma = \frac{1}{K}$), \eqref{discrete} can be rewritten as \[\mathrm{D}_{\mathrm{KL}}(\overbar{\rho}^{n+1} || \overbar{\rho}_\infty) - \mathrm{D}_{\mathrm{KL}}(\overbar{\rho}^n || \overbar{\rho}_\infty) \leq -(1-\gamma)\mathrm{D}_{\mathrm{KL}}(\overbar{\rho}^n || \overbar{\rho}_\infty).\] This can be translated immediately to its continuous-time analog \eqref{entdec}, whence the proof is completed. \qed

\section{Conclusion}
\label{sec:6}
In this manuscript, we have investigated a model (which we call the $K$-averaging model) for a system of self-propelled particles on $\mathbb{R}^d$, in both discrete-time and continuous-time settings. We also provided an rigorous proof on the convergence of the distribution of a typical particle towards a suitable Gaussian equilibrium under the large particle size $N \to \infty$ and large time $n\to \infty$ limit. Even though the majority of the work is done in discrete-time, the relevant results carry over easily to continuous-time. It would also be interesting to examine variants of this model. For instance, the $K$-averaging dynamics on $\mathbb{S}^1$ is closely related to several models in the literature \cite{porfiri_effective_2016,aldana_phase_2003,aldana_phase_2007,pimentel_intrinsic_2008,porfiri_linear_2013}, and it is reasonable to expect a rigorous proof of the corresponding mean-field limit. Unfortunately, the situation on $\mathbb{S}^1$ is inevitably much more complicated since we are lacking the vector-space structure. More generally, averaging is not a straightforward operation over a manifold \cite{degond_local_2014}. Other extensions of the model in the present manuscript are also possible. As of now, every agent communicate with each other. Thus, what would happen if only agents are only interacting through a pre-defined graph of neighboring few chosen neighbors? We would lose the invariance by permutation, thus the notion of limit is more challenging. This would also link the model to certain "consensus models" \cite{hardin_probability_2019,lanchier_probability_2020}. One can also explore different laws of communication between the particles (especially of the non-symmetric and non-all-to-all variety), and investigate the role of noise introduced into the system.\\

\noindent {\bf Acknowledgement} It is a pleasure to thank my Ph.D advisor S\'ebastien Motsch for his tremendous help on various portions of this manuscript.

\bibliographystyle{plain}
\bibliography{K-averaging}

\begin{thebibliography}{10}

\bibitem{aldana_phase_2007}
Maximino Aldana, Victor Dossetti, Christian Huepe, V.~M. Kenkre, and Hernán
  Larralde.
\newblock Phase transitions in systems of self-propelled agents and related
  network models.
\newblock {\em Physical review letters}, 98(9):095702, 2007.
\newblock Publisher: APS.

\bibitem{aldana_phase_2003}
Maximino Aldana and Cristián Huepe.
\newblock Phase transitions in self-driven many-particle systems and related
  non-equilibrium models: a network approach.
\newblock {\em Journal of Statistical Physics}, 112(1-2):135--153, 2003.
\newblock Publisher: Springer.

\bibitem{artstein_solution_2004-1}
Shiri Artstein, Keith Ball, Franck Barthe, and Assaf Naor.
\newblock Solution of {Shannon}’s problem on the monotonicity of entropy.
\newblock {\em Journal of the American Mathematical Society}, 17(4):975--982,
  2004.

\bibitem{bakry_analysis_2013}
Dominique Bakry, Ivan Gentil, and Michel Ledoux.
\newblock {\em Analysis and geometry of {Markov} diffusion operators}, volume
  348.
\newblock Springer Science \& Business Media, 2013.

\bibitem{barbaro_phase_2012}
Alethea~BT Barbaro and Pierre Degond.
\newblock Phase transition and diffusion among socially interacting
  self-propelled agents.
\newblock {\em arXiv preprint arXiv:1207.1926}, 2012.

\bibitem{baumann_modeling_2020}
Fabian Baumann, Philipp Lorenz-Spreen, Igor~M. Sokolov, and Michele Starnini.
\newblock Modeling echo chambers and polarization dynamics in social networks.
\newblock {\em Physical Review Letters}, 124(4):048301, 2020.
\newblock Publisher: APS.

\bibitem{belmonte_self-propelled_2008}
Julio~M. Belmonte, Gilberto~L. Thomas, Leonardo~G. Brunnet, Rita~MC de~Almeida,
  and Hugues Chaté.
\newblock Self-propelled particle model for cell-sorting phenomena.
\newblock {\em Physical Review Letters}, 100(24):248702, 2008.
\newblock Publisher: APS.

\bibitem{berti_almost_2006}
Patrizia Berti, Luca Pratelli, and Pietro Rigo.
\newblock Almost sure weak convergence of random probability measures.
\newblock {\em Stochastics and Stochastics Reports}, 78(2):91--97, 2006.

\bibitem{bertin_boltzmann_2006}
Eric Bertin, Michel Droz, and Guillaume Grégoire.
\newblock Boltzmann and hydrodynamic description for self-propelled particles.
\newblock {\em Physical Review E}, 74(2):022101, 2006.
\newblock Publisher: APS.

\bibitem{bertin_hydrodynamic_2009}
Eric Bertin, Michel Droz, and Guillaume Grégoire.
\newblock Hydrodynamic equations for self-propelled particles: microscopic
  derivation and stability analysis.
\newblock {\em Journal of Physics A: Mathematical and Theoretical},
  42(44):445001, 2009.
\newblock Publisher: IOP Publishing.

\bibitem{billingsley_convergence_2013}
Patrick Billingsley.
\newblock {\em Convergence of probability measures}.
\newblock John Wiley \& Sons, 2013.

\bibitem{boissard_trail_2013}
Emmanuel Boissard, Pierre Degond, and Sebastien Motsch.
\newblock Trail formation based on directed pheromone deposition.
\newblock {\em Journal of mathematical biology}, 66(6):1267--1301, 2013.
\newblock Publisher: Springer.

\bibitem{carlen_kinetic_2013}
Eric Carlen, Robin Chatelin, Pierre Degond, and Bernt Wennberg.
\newblock Kinetic hierarchy and propagation of chaos in biological swarm
  models.
\newblock {\em Physica D: Nonlinear Phenomena}, 260:90--111, 2013.
\newblock Publisher: Elsevier.

\bibitem{chate_modeling_2008}
Hugues Chaté, Francesco Ginelli, Guillaume Grégoire, Fernando Peruani, and
  Franck Raynaud.
\newblock Modeling collective motion: variations on the {Vicsek} model.
\newblock {\em The European Physical Journal B}, 64(3):451--456, 2008.

\bibitem{chuang_state_2007}
Yao-Li Chuang, Maria~R. D’orsogna, Daniel Marthaler, Andrea~L. Bertozzi, and
  Lincoln~S. Chayes.
\newblock State transitions and the continuum limit for a {2D} interacting,
  self-propelled particle system.
\newblock {\em Physica D: Nonlinear Phenomena}, 232(1):33--47, 2007.
\newblock Publisher: Elsevier.

\bibitem{cover_elements_1999}
Thomas~M. Cover.
\newblock {\em Elements of information theory}.
\newblock John Wiley \& Sons, 1999.

\bibitem{dai_pra_stochastic_2017}
Paolo Dai~Pra.
\newblock Stochastic mean-field dynamics and applications to life sciences.
\newblock In {\em International workshop on {Stochastic} {Dynamics} out of
  {Equilibrium}}, pages 3--27. Springer, 2017.

\bibitem{degond_local_2014}
Pierre Degond, Amic Frouvelle, and Gaël Raoul.
\newblock Local stability of perfect alignment for a spatially homogeneous
  kinetic model.
\newblock {\em Journal of Statistical Physics}, 157(1):84--112, 2014.
\newblock Publisher: Springer.

\bibitem{hardin_probability_2019}
Mela Hardin and Nicolas Lanchier.
\newblock Probability of consensus in spatial opinion models with confidence
  threshold.
\newblock {\em arXiv preprint arXiv:1912.06746}, 2019.

\bibitem{hauray_n-particles_2007}
Maxime Hauray and Pierre-Emmanuel Jabin.
\newblock N-particles approximation of the {Vlasov} equations with singular
  potential.
\newblock {\em Archive for rational mechanics and analysis}, 183(3):489--524,
  2007.
\newblock Publisher: Springer.

\bibitem{jabin_quantitative_2018}
Pierre-Emmanuel Jabin and Zhenfu Wang.
\newblock Quantitative estimates of propagation of chaos for stochastic systems
  with {$W^{-1,\infty}$} kernels.
\newblock {\em Inventiones mathematicae}, 214(1):523--591, 2018.
\newblock Publisher: Springer.

\bibitem{lanchier_probability_2020}
Nicolas Lanchier and Hsin-Lun Li.
\newblock Probability of consensus in the multivariate {Deffuant} model on
  finite connected graphs.
\newblock {\em Electronic Communications in Probability}, 25, 2020.
\newblock Publisher: The Institute of Mathematical Statistics and the Bernoulli
  Society.

\bibitem{liggett_interacting_2012}
Thomas~Milton Liggett.
\newblock {\em Interacting particle systems}, volume 276.
\newblock Springer Science \& Business Media, 2012.

\bibitem{madiman_generalized_2007}
Mokshay Madiman and Andrew Barron.
\newblock Generalized entropy power inequalities and monotonicity properties of
  information.
\newblock {\em IEEE Transactions on Information Theory}, 53(7):2317--2329,
  2007.
\newblock Publisher: IEEE.

\bibitem{merle_cutoff_2019}
Mathieu Merle and Justin Salez.
\newblock Cutoff for the mean-field zero-range process.
\newblock {\em The Annals of Probability}, 47(5):3170--3201, 2019.
\newblock Publisher: Institute of Mathematical Statistics.

\bibitem{meleard_propagation_1987}
Sylvie Méléard and Sylvie Roelly-Coppoletta.
\newblock A propagation of chaos result for a system of particles with moderate
  interaction.
\newblock {\em Stochastic processes and their applications}, 26:317--332, 1987.
\newblock Publisher: Elsevier.

\bibitem{motsch_heterophilious_2014}
Sebastien Motsch and Eitan Tadmor.
\newblock Heterophilious dynamics enhances consensus.
\newblock {\em SIAM review}, 56(4):577--621, 2014.

\bibitem{naldi_mathematical_2010}
Giovanni Naldi, Lorenzo Pareschi, and Giuseppe Toscani.
\newblock {\em Mathematical modeling of collective behavior in socio-economic
  and life sciences}.
\newblock Springer Science \& Business Media, 2010.

\bibitem{oelschlager_martingale_1984}
Karl Oelschlager.
\newblock A martingale approach to the law of large numbers for weakly
  interacting stochastic processes.
\newblock {\em The Annals of Probability}, pages 458--479, 1984.
\newblock Publisher: JSTOR.

\bibitem{pimentel_intrinsic_2008}
Jaime~A. Pimentel, Maximino Aldana, Cristián Huepe, and Hernán Larralde.
\newblock Intrinsic and extrinsic noise effects on phase transitions of network
  models with applications to swarming systems.
\newblock {\em Physical Review E}, 77(6):061138, 2008.
\newblock Publisher: APS.

\bibitem{popoviciu_sur_1935}
Tiberiu Popoviciu.
\newblock Sur les équations algébriques ayant toutes leurs racines réelles.
\newblock {\em Mathematica}, 9:129--145, 1935.

\bibitem{porfiri_linear_2013}
Maurizio Porfiri.
\newblock Linear analysis of the vectorial network model.
\newblock {\em IEEE Transactions on Circuits and Systems II: Express Briefs},
  61(1):44--48, 2013.
\newblock Publisher: IEEE.

\bibitem{porfiri_effective_2016}
Maurizio Porfiri and Gil Ariel.
\newblock On effective temperature in network models of collective behavior.
\newblock {\em Chaos: An Interdisciplinary Journal of Nonlinear Science},
  26(4):043109, 2016.

\bibitem{rezakhanlou_entropy_2008}
Fraydoun Rezakhanlou, Cédric Villani, and François Golse.
\newblock {\em Entropy methods for the {Boltzmann} equation: lectures from a
  special semester at the {Centre} Émile {Borel}, {Institut} {H}. {Poincaré},
  {Paris}, 2001}.
\newblock Number 1916. Springer Science \& Business Media, 2008.

\bibitem{sznitman_topics_1991}
Alain-Sol Sznitman.
\newblock Topics in propagation of chaos.
\newblock In {\em Ecole d'été de probabilités de {Saint}-{Flour}
  {XIX}—1989}, pages 165--251. Springer, 1991.

\end{thebibliography}

\end{document}